\documentclass[11pt]{amsart}
\usepackage{amssymb}
\usepackage{eucal}

\newcommand{\eps}{\varepsilon}
\newcommand\cc{C_{\alpha\beta\gamma}^{ijk}}

\newcommand\od{\kappa}
\newcommand\rr{{\mathbb R}}

\newtheorem{lemma}{Lemma}[section]
\newtheorem{theorem}{Theorem}[section]

\theoremstyle{remark}
\newtheorem*{remark}{Remark}

\theoremstyle{definition}

\author{Thomas C. Sideris}
\address{Department of Mathematics\\
University of California\\
Santa Barbara, CA 93106}
\email{sideris@math.ucsb.edu}
\thanks{T.C.S. was supported in part by the National Science Foundation.}

\author{Shu-Yi Tu}
\address{Department of Mathematics\\
St.\ Cloud State University\\
St.\ Cloud, MN 56301}
\email{stu@stcloudstate.edu}
\title[Nonlinear Wave Equations]
{Global Existence for Systems of Nonlinear Wave
Equations in 3D with Multiple Speeds}

\begin{document}

\begin{abstract}{Global smooth solutions to the initial value
problem for systems of nonlinear
wave equations with multiple propagation speeds will
be constructed in the case of small initial data and nonlinearities
satisfying the null condition.}
\end{abstract}



\maketitle
\numberwithin{equation}{section}

\section{Introduction}

This paper is concerned with the Cauchy problem for 
coupled systems of quasilinear wave equations in three space
dimensions of the form
\[
\partial_t^2 u^k-c_k^2\triangle u^k
=C^{jk}_{\alpha\beta}(\partial u)
\partial_\alpha\partial_\beta u^j,\quad k=1,\ldots,m,
\]
subject to suitably small initial conditions.
We assume that the propagation speeds are distinct, 
and we refer to this situation as the nonrelativistic case.
Here, $\partial u$ stands for the full space-time gradient,
and $C^{jk}_{\alpha\beta}(\xi)=
{O}(|\xi|)$ are
smooth functions near the
origin in $\rr^{4m}$.
We shall construct a unique global classical solution,
provided that the coefficients of the nonlinear terms
satisfy the null condition.
This nonrelativistic system serves as a simplified model
for wave propagation problems with different speeds,
such as nonlinear elasticity, charged plasmas,
and magneto-hydrodynamics.

The main difficulty in the nonrelativistic case is
that the smaller symmetry group of the linear operator
weakens the form of the invariant Klainerman
inequality, see Section \ref{sob}.  In order to obtain a viable
$L^\infty-L^2$ estimate for solutions,
we utilize an additional set of weighted $L^2$
estimates, as has been developed in \cite{KlSid},
\cite{Sid}, \cite{Sid2}.  The advantage of this
method is the total avoidance of direct estimation
of the fundamental solution for the linear problem
or any type of asymptotic constructions.
We treat nondivergence form nonlinearities which
may contain both spatial and temporal derivatives.

In the 3D relativistic (scalar) case, the null condition was
first identified and shown to lead to global existence
of small solutions in \cite{Ch}, \cite{Kl2}.
Without it, small solutions remain smooth ``almost
globally" \cite{Kl1}, but arbitrarily small initial
conditions can develop singularities in finite time
\cite{J}, \cite{Sid1}.  Small solution always exist
globally in higher dimensions \cite{KlPon}, \cite{Sh}, \cite{Kl1}.
The 2D relativistic case is rather more complicated.
The sharpest results are given in \cite{Al},
but other work appeared previously in \cite{Ho}, \cite{Ka}.

The case of nonrelativistic systems has been considered
in  3D \cite{yo}  and in 2D \cite{Kubo}, \cite{Ho1}.
We mention also the early work \cite{Kov}, \cite{Kov1}, \cite{Kov2}
which deals with nonresonant interactions.
The common theme in these works is the direct estimation
of the fundamental solution which, as mentioned above, is
avoided here.  

The statement of the main result is given in section \ref{mainres} 
after the introduction of some standard notation.  
The rest of the paper presents the proof.
To simplify the exposition, we truncate the nonlinearity at the
quadratic level, but this entails no loss of generality
since the higher-order terms do not affect the global
behavior of small solutions, \cite{Kl1}.


\bigskip
\section{Notation}

Points in $\rr^4$ will be denoted by $X=(x^0,x^1,x^2,x^3)=(t,x)$.
Partial derivatives will be written as
$ \partial_k={\partial}/{\partial x^k}$,
$k=0,\ldots,3 $, 
with the abbreviations  
$ \partial=(\partial_0,\partial_1,\partial_2,\partial_3)=(\partial_t,\nabla)$.
The angular-momentum operators are defined as
\begin{equation*}
\Omega=(\Omega_1,\Omega_2,\Omega_3)=x\wedge\nabla,
\end{equation*}
where $\wedge$ denotes the usual vector cross product in $\rr^3$, and
the scaling operator is defined by
\begin{equation}
\label{scaling}
S=t\partial_t+r\partial_r=x^\alpha\partial_\alpha.
\end{equation}
The collection of these seven vector fields will be labeled as
\[
\Gamma=(\Gamma_0,\ldots,\Gamma_7)=(\partial,\Omega,S).
\]
Instead of the usual multi-index notation, we will write
$a=(a_1,\ldots,a_\od)$ for a sequence of indices $a_i\in
\{0,\ldots,7\}$ of length $|a|=\od$, and 
\[
\Gamma^a=\Gamma_{a_\od}
\cdots\Gamma_{a_1}.
\]
Suppose that $b$ and $c$ are disjoint subsequences of $a$.
Then we will say $b+c=a$, if $|b|+|c|=|a|$,  and $b+c<a$,
if $|b|+|c|<|a|$.

The d'Alembertian will be used to denote the operator
\[
\square=\mbox{Diag}(\square_1,\ldots,\square_m)
\quad\mbox{with}\quad
\square_k=\partial^2_t-c^2_k\triangle.
\]
For convenience, we will assume that the speeds are distinct
\[
c_1>\cdots>c_m>0.
\] 
It is also possible to treat the case where some of the
speeds are the same, see the remark following the statement
of Theorem \ref{thm}.
This operator acts on vector functions $u:\rr^4\to\rr^m$.
The standard energy  then is defined as
\[E_1(u(t))=\sum_{k=1}^{m}\int_{{\rr}^3}
[\,|\partial_tu^k(t,x)|^2+c_k^2\,|\nabla u^k(t,x)|^2\,]\,dx,
\]
and higher order derivatives will be estimated through
\begin{subequations}
\begin{equation}
\label{ennorm}
E_\od(u(t))=\sum_{|a|\le\od-1} E_1(\Gamma^au(t)),
\qquad \od=2, 3, \ldots
\end{equation}

In order to describe the solution space, we introduce the
time-independent vector fields
$\Lambda=(\Lambda_1,\ldots,\Lambda_7)=(\nabla,\Omega,r\partial_r)$.
Define 
\[
H^\od_\Lambda(\rr^3)=\{f\in L^2(\rr^3;\rr^m) : \Lambda^af\in L^2,
\; |a|\le\od\},
\]
with the norm
\begin{equation}
\label{sobnorm}
\|f\|_{H^\od_\Lambda}=\sum_{|a|\le\od}\|\Lambda^af\|_{L^2}.
\end{equation}

Solutions will be constructed in the space
$\dot{H}_\Gamma^\od(T)$ obtained by closing
the set $C^\infty([0,T);C_0^\infty(\rr^3,\rr^m))$
in the norm $\sup\limits_{0\le t< T} E_\od^{1/2}(u(t))$.
Thus,
$$
\dot{H}_\Gamma^\od(T)\subset\left\{u(t,x) : \partial u(t,\cdot)
\in\bigcap_{j=0}^{\od-1} C^j([0,T);H^{\od-1-j}_\Lambda)\right\}.
$$
By \eqref{ellinf4}, it will follow that
$\dot{H}_\Gamma^\kappa(T)\subset C^{\kappa-2}([0,T)\times\rr^3;\rr^m)$.

An important intermediate role will played by the weighted norm
\begin{equation}
\label{wnorm}
\mathcal{X}_\od(u(t))=\sum_{k=1}^m\sum_{|a|=2}\sum_{|b|\le\od-2}
\|\langle c_kt-|x|\rangle\partial^a\Gamma^bu^k(t)\|_{L^2(\rr^3)},
\end{equation}
\end{subequations}
where we use the notation $\langle\rho\rangle=(1+|\rho|^2)^{1/2}$.

\section{Main Result}
\label{mainres}

Consider the initial value problem
for a coupled nonlinear system of the form
\begin{equation}
\label{pde}
\square u = N(u,u)
\end{equation}
in which the components of the quadratic nonlinearity depend on
the form
\begin{subequations}
\begin{equation}
\label{non}
N^k(u,v)=\cc \partial_\alpha u^i\partial_\beta\partial_\gamma v^j.
\end{equation}
Summation is performed over repeated indices regardless of their
position, up or down.  Greek indices range from 0 to 3 and Latin
indices from 1 to $m$.  

Existence of solutions depends on the energy method which
requires the system to be symmetric:
\begin{equation}
\label{symmetric}
\cc=C^{ikj}_{\alpha\beta\gamma}=C^{ijk}_{\alpha\gamma\beta}.
\end{equation}

The key assumption necessary for global existence is the following
{\em null} condition which says that the self-interaction 
of each wave family is nonresonant:
\begin{equation}
\label{nc}
C^{kkk}_{\alpha\beta\gamma}X_\alpha X_\beta X_\gamma=0
\quad\mbox{for all}\quad X\in \mathcal{N}_k,
\quad k=1,\ldots,m,
\end{equation}
\end{subequations}
with the null cones
\[
\mathcal{N}_k=\{X\in\rr^4: X_0^2-c_k^2(X_1^2+X_2^2+X_3^2)=0\}.
\]

\begin{theorem}
\label{thm}
Assume that the nonlinear terms in \eqref{non}
satisfy the symmetry and null conditions \eqref{symmetric}, \eqref{nc}.
Then the initial value problem for \eqref{pde}
with initial data
\[
\partial_\alpha u(0)\in H^{\od-1}_\Lambda(\rr^3),\quad \od\ge 9
\]
satisfying 
\begin{equation}
\label{smallness}
E_{\od-2}^{1/2}(u(0))\;\exp\;CE_\od^{1/2}(u(0)) < \eps,
\end{equation}
with $\eps$ sufficiently small,
has a unique global solution $u\in \dot H_\Gamma^\od(T)$ for
every $T>0$. 
The solution satisfies the bounds
\[
E_{\od-2}^{1/2}(u(t))< 2\eps\quad\mbox{and}
\quad E_\od(u(t))\le 4 E_\od(u(0)) \langle t\rangle^{C{\eps}}.
\]
\end{theorem}

\begin{remark}
We briefly discuss the case when some of the speeds are repeated.
Suppose that only $\ell<m$ of the speeds, $c_1>c_{k_2}>\ldots
>c_{k_\ell}$ are distinct.  For $p=1,\ldots,\ell$,
let $I_p=\{k:1\le k\le m,\;  c_k=c_{k_p}\}.$
The null condition is now extended to be
\[
C^{ijk}_{\alpha\beta\gamma}X_\alpha X_\beta X_\gamma=0,
\quad \mbox{for all}\quad X\in \mathcal{N}_{k_p},\;
(i,j,k)\in I_p^3,\; p=1,\ldots,\ell.
\]
The proof can easily be adjusted to handle this more
general case.
\end{remark}


\section{Commutation and Null Forms}

In preparation for the energy estimates, we need to consider
the commutation properties of the vector fields $\Gamma$ with respect
to the nonlinear terms.  It is necessary to verify that the
null structure is preserved upon differentiation.


\begin{lemma}
\label{nd}
Let $u$ be solution $u$ of \eqref{pde} in $\dot{H}^\od_\Gamma(T)$.
Assume that the null condition \eqref{nc} holds for the nonlinearity
in \eqref{non}.
Then for $|a|\le\od-1$,
\[
\square \Gamma^a u = \sum_{b+c+d=a}N_d(\Gamma^bu,\Gamma^cu),
\]
in which each $N_d$ is a quadratic nonlinearity of the form
\eqref{non} satisfying \eqref{nc}.
Moreover, if $b+c=a$, then $N_d=N$.
\end{lemma}

\begin{proof}
First we note the well-known facts that
\[
[\partial,\square]=0,\quad[\Omega,\square]=0,\quad[S,\square]=-2\square.
\]
Recalling the definition \eqref{non}, we set
\[
[\Gamma,N](u,v)=\Gamma N(u,v) - N(\Gamma u,v)-N(u,\Gamma v).
\]
This is a quadratic nonlinearity of the form \eqref{non}.
Thus, if $[\Gamma,N]$ is null for each $\Gamma$, then the result follows
by induction.
In fact, if $d=(d_1,\ldots,d_k)$, then $N_d$ is the $k$-fold
commutator $N_d=[\Gamma_{d_k},[\ldots,[\Gamma_{d_1},N]]]$.

A simple calculation shows that
\[
[\partial,N](u,v)=0\quad\mbox{and}\quad[S,N]=-3N(u,v).
\]
Thus, these commutators are null if $N$ is null.

We can express the angular momentum operators as
$\Omega_\lambda=\eps_{\lambda\mu\nu}x_\mu\partial_\nu$,
$\lambda=1,2,3$, where
$\eps_{\lambda\mu\nu}$ is the tensor with
value $+1$, $-1$ if $\lambda\mu\nu$ is an even, respectively odd,
 permutation of 123, and with value 0 otherwise.  
Using this, we find that the $k^{\mbox{\em th}}$ component of 
$[\Omega_\lambda,N]$ is
\[
[\Omega_\lambda,N]^k(u,v)
=\widetilde{C}^{ijk}_{\alpha\beta\gamma}\partial_\alpha u^j
\partial_\beta\partial_\gamma v^k
\]
with
\[
\widetilde{C}^{ijk}_{\alpha\beta\gamma}
=[C^{ijk}_{\alpha\beta\nu}\eps_{\lambda\gamma\nu}
+C^{ijk}_{\nu\beta\gamma}\eps_{\lambda\alpha\nu}
+C^{ijk}_{\alpha\nu\gamma}\eps_{\lambda\beta\nu}].
\]
To see that this commutator is also null, write
\[
h^k(X)=C^{kkk}_{\alpha\beta\gamma}X_\alpha X_\beta X_\gamma
\quad\mbox{and}\quad
\tilde{h}^k(X)=\widetilde{C}^{kkk}_{\alpha\beta\gamma}X_\alpha
X_\beta X_\gamma.
\]
Then $\tilde{h}^k(X)=Dh^k(X)Y^\lambda$ with
$Y^{\lambda}_\mu=\eps_{\lambda\mu\nu}X_\nu$.
Now the null condition says that $h^k(X)=0$ for $X\in\mathcal{N}_k$.
But since $Y^\lambda$ is tangent to $\mathcal{N}_k$ at $X$,  we have
$\tilde{h}^k(X)=0$ for $X\in\mathcal{N}_k$. This implies that
$[\Omega_\lambda,N]$ is null.
\end{proof}

\section{Estimates for Null Forms}

The utility of the null condition is captured in the next lemma.
The presence of the the terms with the weight $\langle c_kt-r\rangle$
in these inequalities is explained by the absence of the
Lorentz rotations in our list of vector fields $\Gamma$.

\begin{lemma}
\label{ptwnc}
Suppose that the nonlinear form $N(u,v)$
defined \eqref{non} satisfies the null condition \eqref{nc}.  
Set $c_0=\min\{c_k/2: k=1,\ldots,m\}$.  For
$u$, $v$, $w\in C^2([0,T]\times\rr^3;\rr^m)$ and $r\ge c_0 t$, we have
at any point $X=(t,x)$
\begin{subequations}
\begin{multline}
\label{ptwnc1}
|C^{kkk}_{\alpha\beta\gamma}\partial_\alpha u^k 
\partial_\beta\partial_\gamma v^k|\\
\le\frac{C}{\langle X\rangle}\Big[|\Gamma u^k||\partial^2v^k|
+|\partial u^k||\partial\Gamma v^k|
+\langle c_kt-r\rangle|\partial u^k|
|\partial^2v^k| \Big]
\end{multline}
and
\begin{multline}
\label{ptwnc2}
|C^{kkk}_{\alpha\beta\gamma}\partial_\alpha u^k
\partial_\beta v^k \partial_\gamma w^k|\\
\le \frac{C}{\langle X\rangle}\Big[
|\Gamma u^k| |\partial v^k| |\partial w^k|
+ |\partial u^k| |\Gamma v^k| |\partial w^k|
+ |\partial u^k| |\partial v^k| |\Gamma w^k|\\
+ \langle c_kt-r\rangle |\partial u^k| |\partial v^k| |\partial w^k|\Big],
\end{multline}
\end{subequations}
in which $ \langle X\rangle = (1+|X|^2)^{1/2}$.
\end{lemma}
\begin{proof}

Spatial derivatives have the decomposition
\[
\nabla=\frac{x}{r}\partial_r - \frac{x}{r^2}\wedge\Omega.
\]
So if we introduce the two operators 
$D_k^{\pm}=\frac{1}{2}\,(\partial_t\pm c_k \partial_r)$
and the null vectors $Y^\pm_k=(1,\pm x/c_kr)\in\mathcal{N}_k$
we obtain
\begin{equation}
\label{dd1}
(\partial_t,\nabla)=(Y^-_kD^-_k+Y^+_kD^+_k)
-\left(0,\frac{x}{c_kr^2}\wedge\Omega\right).
\end{equation}
On the other hand, if we write
\[
D_k^+=\frac{c_k}{c_kt+r}S - \frac{c_kt-r}{c_kt+r}D_k^-,
\]
the formula \eqref{dd1} can be transformed into
\[
\partial= Y_k^-D_k^--\frac{c_kt-r}{c_kt+r}Y_k^+D_k^- +\frac{c_k}{c_kt+r}Y_k^+S
-\left(0,\frac{x}{c_kr^2}\wedge\Omega\right).
\]
Thus, we have
\begin{subequations}
\begin{equation}
\label{dd21}
\partial \equiv Y_k^- D_k^- + R.
\end{equation}
Now, we may assume that $|X|\ge1$, for otherwise the estimates
are trivial.  But then it follows that $1/r$ and
$1/(c_kt+r)$ are bounded by $C/\langle X\rangle$,
and as a consequence we have
\begin{equation}
\label{dd22}
|Ru|\le C\langle X\rangle^{-1} [ |\Gamma u| + \langle c_kt-r\rangle
|\partial u|]
\end{equation}
\end{subequations}

Using \eqref{dd21}, we have
\begin{multline}
\label{expand}
C^{kkk}_{\alpha\beta\gamma}
\partial_\alpha u^k \partial_\beta\partial_\gamma v^k
=C^{kkk}_{\alpha\beta\gamma}[Y^-_{k\alpha} Y^-_{k\beta} Y^-_{k\gamma}
D_k^-u^k(D_k^-)^2v^k\\
+R_\alpha u^k\partial_\beta\partial_\gamma v^k
+Y^-_{k\alpha}D_k^-u^k R_\beta\partial_\gamma v^k
+Y^-_{k\alpha} D_k^-u^kY^-_{k\beta}D_k^-R_\gamma v^k].
\end{multline}
The first term in \eqref{expand}
vanishes since $N$ obeys the null condition, and by \eqref{dd22}
the remaining terms in \eqref{expand} have the estimate
\eqref{ptwnc1}.

The proof of \eqref{ptwnc2} is similar.

\end{proof}

\section{Sobolev Inequalities}
\label{sob}

The following Sobolev inequalities involve only
the angular momentum operators since we are
in the nonrelativistic case.  The weight $\langle ct-r\rangle$
compensates for this.  We use the notation defined in 
\eqref{ennorm}, \eqref{sobnorm}, \eqref{wnorm}.

\begin{lemma}
Let $u\in \dot{H}_\Gamma^\od(T)$, with $\mathcal{X}_\od(u(t))
<\infty$.
\begin{alignat}{2}
\label{ellinf4}
&\langle r \rangle^{1/2}|\Gamma^a u(t,x)|\le CE_\od^{1/2}(u(t)),
&\quad &|a|+2\le\od\\
\label{ellinf5}
&\langle r \rangle |\partial\Gamma^a u(t,x)|\le CE_\od^{1/2}(u(t)),
&\quad &|a|+3\le\od\\
\label{ellinf6}
&\langle r\rangle\langle c t -r \rangle ^{1/2}
|\partial\Gamma^a u(t,x)|\\
\nonumber
&\qquad\le C\Big[E_\od^{1/2}(u(t))+\mathcal{X}_\od(u(t))\Big],
&\quad &|a|+3\le\od\\
\label{ellinf7}
&\langle r\rangle\langle c t -r \rangle
| \partial^2\Gamma^a u(t,x)|
\le C\mathcal{X}_\od(u(t)),
&\quad&|a|+4\le\od.
\end{alignat}
\end{lemma}

\begin{proof}
This  result is essentially Proposition 3.3 in  \cite{Sid2}.
\end{proof}

\section{Weighted Decay Estimates}

The main extra step in the nonrelativistic case
is to control the weighted norm $\mathcal{X}_\od(u(t))$.
This will be accomplished in this section by a type of
bootstrap argument.  
\begin{lemma}
\label{decay1}
Let $u\in \dot{H}_\Gamma^\od(T)$.  Then 
\begin{equation}
\label{locen}
\mathcal{X}_\od(u(t))\le C\left[E_{\od}^{1/2}(u(t))+
\sum_{|a|\le\od-2}\|(t+r)\,\square\Gamma^a u(t)\|^2_{L^2}\right].
\end{equation}
\end{lemma}

\begin{proof}
Recall that the weighted norm involves derivatives
in the form $\partial^2\Gamma^a u$
In the case when $\partial^2=\nabla\partial$,
the result was given in Lemma 3.1 of \cite{KlSid}.  Otherwise,
if $\partial^2=\partial_t^2$,
then the result is an immediate consequence of (2.10) in \cite{KlSid}.
\end{proof}

Now we assume that $u$ solves the nonlinear PDE.

\begin{lemma}
\label{lem2}
Let $u\in \dot{H}_\Gamma^\od(T)$ be a solution of \eqref{pde}.
Define ${\od}'=\left[\frac{\od-1}{2}\right]+3.$  
Then  for all $|a|\le\od-2$,
\begin{multline}
\label{w2}
\|(t+r)\square\Gamma^a u(t)\|^2_{L^2}\\
\le C[\mathcal{X}_{\od'}(u(t))
E_{\od}^{1/2}(u(t))
+\mathcal{X}_{\od}(u(t))
E_{{\od}'}^{1/2}(u(t))].
\end{multline}
\end{lemma}
\begin{proof}
By Lemma \ref{nd}, we must estimate terms of the form
\[
\|(t+r)\partial\Gamma^bu^i\partial^2\Gamma^cu^j\|_{L^2},
\]
but since $(t+r)\le C\langle r \rangle \langle c_jt-r\rangle$,
we will consider
\begin{equation}
\label{expression}
\|\langle r \rangle\langle c_jt-r\rangle
\partial\Gamma^bu^i
\partial^2\Gamma^cu^j\|_{L^2},
\end{equation}
with $b+c\le a$, and $|a|\le\od-2$.  

Let $m=\left[\frac{\od-1}{2}\right]
=\od'-3$.
We separate two cases:  either $|b|\le m$ or $|c|\le m-1$.
In the first case, \eqref{expression} is estimated as follows
using \eqref{ellinf5}:
\[
\|\langle r \rangle\partial\Gamma^bu^i\|_{L^\infty}
\|\langle c_jt-r\rangle\partial^2\Gamma^cu^j\|_{L^2}
\le CE_{\od'}^{1/2}(u(t))\mathcal{X}_\od(u(t)).
\]

Otherwise, we use \eqref{ellinf7} to estimate \eqref{expression} by:
\[
\|\partial\Gamma^bu^i\|_{L^2}
\|\langle r \rangle\langle c_jt-r\rangle\partial^2\Gamma^cu^j\|_{L^\infty}
\le C E_{\od}^{1/2}(u(t))\mathcal{X}_{\od'}(u(t)).
\]
\end{proof}

The next result gains control of the weighted norm by the energy.
We distinguish two different energies, the smaller of which must
remain small.  In the next section, we will
allow the larger energy will grow polynomially in time.

\begin{lemma}
\label{we}
Let $u\in \dot{H}_\Gamma^\od(T)$, $\od\ge8$, be a solution of \eqref{pde}.
Define $\mu=\od-2$, and assume that 
\[
\eps_0\equiv\sup_{0\le t < T} E_{\mu}^{1/2}(u(t))
\]
is sufficiently small.
Then for $0\le t< T$, 
\begin{subequations}
\begin{equation}
\label{we1}
\mathcal{X}_\mu(u(t))\le C E_{\mu}^{1/2}(u(t))
\end{equation}
and
\begin{equation}
\label{we2}
\mathcal{X}_\od(u(t))\le C E_{\od}^{1/2}(u(t)).
\end{equation}
\end{subequations}
\end{lemma}
\begin{proof}
Let $\mu' = \left[\frac{\mu-1}{2}\right]+3$, $\mu=\od-2$.
Since $\mu\ge 6$, we have $\mu'\le\mu$.  Thus, by Lemmas
\ref{decay1} and \ref{lem2}, we find using our assumption
\[
\mathcal{X}_\mu(u(t))\le C [E_{\mu}^{1/2}(u(t)) +  \eps_0 
\mathcal{X}_\mu(u(t))].
\]
Thus, if $\eps_0$ is small enough, the bound \eqref{we1} results.

Again since $\od\ge8$, we have $\od'=\left[\frac{\od-1}{2}\right]+3
\le \mu =\od-2$.   From Lemmas \ref{decay1} and \ref{lem2} we now
have
\[
\mathcal{X}_\od(u(t))\le C [E_{\od}^{1/2}(u(t))
+\mathcal{X}_{\mu}(u(t))
E_{\od}^{1/2}(u(t))
+\mathcal{X}_{\od}(u(t))
E_{\mu}^{1/2}(u(t))].
\]
If we apply \eqref{we1} and our assumption, then
\[
\mathcal{X}_\od(u(t))\le C [E_{\od}^{1/2}(u(t))
+\eps_0\mathcal{X}_{\od}(u(t))],
\]
from which \eqref{we2} follows.
\end{proof}

\section{Energy Estimates}

\subsection*{General energy method}
In this section we shall complete the proof of Theorem \ref{thm}.
Assume that $u(t)\in \dot{H}^\od_\Gamma(T)$ is a local
solution of the initial value problem for \eqref{pde}.
Our task will be to show that $E_\od(u(t))$ remains finite
for all $t\ge0$.  To do so, we will derive a pair
of coupled differential inequalities for (modifications of)
$E_\od(u(t))$ and $E_{\mu}(u(t))$, with $\mu=\od-2$.
If \eqref{smallness} holds
then $E_\mu^{1/2}(u(0))<\eps$. 
Suppose that $T_0$ is the largest
time such that  $E_\mu^{1/2}(u(t))<2\eps$,
for $0\le t< T_0$ with $\eps$  small enough so
that Lemma \ref{we} is valid.
All of the following computations will be valid on
this time interval.

Following the energy method, we have for any $\nu=1,\ldots,\od$,
\[
E_\nu'(u(t))=\sum_{|a|\le\nu-1}\int\langle\square\Gamma^au(t),
\partial_t\Gamma^au(t)\rangle dx,
\]
and from Lemma \ref{nd}, this takes the form
\begin{equation}
\label{enmain}
E_\nu'(u(t))=\sum_{|a|\le\nu-1}\sum_{\phantom{|}b+c+d=a}
\int \langle N_d(\Gamma^bu,\Gamma^cu),\partial_t\Gamma^au\rangle dx.
\end{equation}

Terms in \eqref{enmain} with $b=0$, $c=a$, and $|a|=\nu-1$ are handled
with the aid of the symmetry condition \eqref{symmetric}
which allows us to integrate by parts as follows.
Recall that from Lemma \ref{nd}, $N_d=N$ when $b+c=a$.
\begin{align*}
\int \langle N(u,\Gamma^au),\partial_t\Gamma^au\rangle dx
=& \cc \int \partial_\alpha u^i\partial_\beta\partial_\gamma u^j
\partial_t u^kdx\\
=& \cc \int \partial_\gamma[\partial_\alpha u^i\partial_\beta \Gamma^au^j
\partial_t \Gamma^au^k]dx\\
&-\cc \int \partial_\alpha\partial_\gamma u^i\partial_\beta \Gamma^au^j
\partial_t \Gamma^au^kdx\\
&-\cc \int \partial_\alpha u^i\partial_\beta \Gamma^au^j
\partial_t \partial_\gamma \Gamma^au^kdx\\
=& C^{ijk}_{\alpha\beta 0} 
\partial_t\int \partial_\alpha u^i\partial_\beta \Gamma^au^j
\partial_t \Gamma^au^kdx\\
&-\cc \int \partial_\alpha\partial_\gamma u^i\partial_\beta \Gamma^au^j
\partial_t \Gamma^au^kdx\\
&-\frac{1}{2}\;\cc 
\int \partial_\alpha u^i\partial_t[\partial_\beta \Gamma^au^j
\partial_\gamma \Gamma^au^k] dx\\
=&\frac{1}{2}\;\cc\eta_{\gamma\delta}\;
\partial_t\int \partial_\alpha u^i\partial_\beta \Gamma^au^j
\partial_\delta \Gamma^au^kdx\\
&-\cc \int \partial_\alpha\partial_\gamma u^i\partial_\beta \Gamma^au^j
\partial_t \Gamma^au^kdx\\
&+\frac{1}{2}\;\cc\int\partial_t\partial_\alpha u^i\partial_\beta \Gamma^au^j
\partial_\gamma \Gamma^au^kdx,
\end{align*}
using the symbol $\eta_{\gamma\delta}=\mbox{Diag}[1,-1,-1,-1]$.
The first term above can be absorbed into the energy as a lower order
perturbation.  Define
\[
\widetilde{E}_\nu(u(t))=E_\nu(u(t))
- \frac{1}{2}\sum_{|a|=\nu-1}\cc\eta_{\gamma\delta}\;
\int \partial_\alpha u^i\partial_\beta \Gamma^au^j
\partial_\delta \Gamma^au^kdx.
\]
The perturbation is bounded by $C\|\nabla u\|_{L^\infty}E_\nu(u(t))$,
but by \eqref{ellinf5}, the maximum norm
$\|\nabla u\|_{L^\infty}$ is controlled
by $E_3^{1/2}\le E_\mu^{1/2}<2\eps$.
Thus, for small solutions we have
\begin{equation}
\label{enequiv}
(1/2)E_\nu(u(t))\le \widetilde{E}_\nu(u(t)) \le 2E_\nu(u(t)).
\end{equation}

Returning to \eqref{enmain}, we have derived
the energy identity
\begin{align}
\label{enmain1}
\widetilde{E}'_\nu(u(t))=&
\sum_{|a|\le\nu-1\phantom{\begin{smallmatrix}|\\1\end{smallmatrix}}}
\sum_{\begin{smallmatrix}b+c+d=a\\|a|\ne\nu-1\end{smallmatrix}}
\int \langle N_d(\Gamma^bu,\Gamma^cu),\partial_t\Gamma^au\rangle dx\\
\nonumber
&
+\sum_{|a|=\nu-1\phantom{\begin{smallmatrix}|\\1\end{smallmatrix}}}
\left[\sum_{\begin{smallmatrix}b+c=a\\c\ne a\end{smallmatrix}}\int \langle
N(\Gamma^bu,\Gamma^cu),\partial_t\Gamma^au\rangle dx\right.\\
\nonumber
&
\phantom{+\sum_{b+c+d=a}\Bigg[}
-\cc \int \partial_\alpha\partial_\gamma u^i\partial_\beta \Gamma^au^j
\partial_t \Gamma^au^kdx\\
\nonumber
&
\left.
\phantom{+\sum_{b+c+d=a}\Bigg[}
+\frac{1}{2}\;\cc\int\partial_t\partial_\alpha u^i\partial_\beta
\Gamma^au^j
\partial_\gamma \Gamma^au^kdx\right].
\end{align}

\subsection*{Higher energy}
For the first series of estimates we take $\nu=\od$ in \eqref{enmain1}.
We obtain immediately
\begin{equation}
\label{en1}
\widetilde{E}'_\od(u(t))\le C 
\sum_{\phantom{i}i,j,k\phantom{\begin{smallmatrix}i\\|\end{smallmatrix}}}
\sum_{\phantom{i}|a|\le\od-1\phantom{\begin{smallmatrix}i\\|\end{smallmatrix}}}
\sum_{\begin{smallmatrix}b+c\le a\\ c\ne a\end{smallmatrix}}
\|\partial\Gamma^b u^i \partial^2\Gamma^cu^j\|_{L^2}
\|\partial\Gamma^au^k\|_{L^2}.
\end{equation}
In some cases, the indices $i$ and $j$ have been interchanged.
In the sum on the right-hand side of \eqref{en1},
 we have either $|b|\le \od'$
or $|c|\le \od'-1$, with $\od'=\left[\frac{\od}{2}\right]$.
Note that since $\od\ge 9$, we have $\od'+3\le\od-2=\mu$.
We will also use that $\langle t \rangle 
\le C \langle r \rangle \langle c_jt-r \rangle$.

In the first case, we estimate using \eqref{ellinf5} and \ref{we2}
\begin{align*}
\|\partial\Gamma^b u^i \partial^2\Gamma^cu^j\|_{L^2}
\le & C\langle t \rangle^{-1} 
\|\langle r \rangle \partial\Gamma^b u^i\|_{L^\infty}
\|\langle c_jt-r \rangle \partial^2\Gamma^cu^j\|_{L^2}\\
\le & C\langle t \rangle^{-1} E_{|b|+3}^{1/2}(u(t)) \mathcal{X}_\od(u(t))\\
\le & C\langle t \rangle^{-1} E_\mu^{1/2}(u(t)) E_\od^{1/2}(u(t)).
\end{align*}

In the second case, we use \eqref{ellinf7} and then \eqref{we1}
\begin{align*}
\|\partial\Gamma^b u^i \partial^2\Gamma^cu^j\|_{L^2}
\le & C\langle t \rangle^{-1}
\|\partial\Gamma^b u^i\|_{L^2}
\|\langle r \rangle\langle c_jt-r \rangle \partial^2\Gamma^cu^j\|_{L^\infty}\\
\le & C\langle t \rangle^{-1}
E_\od^{1/2}(u(t)) \mathcal{X}_{|c|+4}(u(t))\\
\le & C\langle t \rangle^{-1}
E_\od^{1/2}(u(t)) \mathcal{X}_{\mu}(u(t))\\
\le & C\langle t \rangle^{-1}
E_\od^{1/2}(u(t)) E_\mu^{1/2}(u(t)).
\end{align*}

Going back to \eqref{en1} and recalling \eqref{enequiv},
we have established the inequality
\begin{align}
\label{en2}
\widetilde{E}'_\od(u(t))
&\le C \langle t \rangle^{-1} E_\mu^{1/2}(u(t)) E_\od(u(t))\\
\nonumber
&\le C \langle t \rangle^{-1} \widetilde{E}_\mu^{1/2}(u(t)) 
\widetilde{E}_\od(u(t)).
\end{align}

\subsection*{Lower energy}
The second series of energy estimates will exploit the
null condition.   We return  to  \eqref{enmain1} now with
$\nu=\mu=\od-2$.
The resulting integrals on the right-hand
side of \eqref{enmain1}
will be subdivided into separate integrals over the
regions $r\le c_0 t$ and $r\ge c_0t$.  Recall that the constant
$c_0$ was defined in Lemma \ref{ptwnc}.

\subsubsection*{Inside the cones}
On the region $r\le c_0 t$, we have that the right-hand side
of \eqref{enmain1} is bounded above by
\[
\sum_{i,j,k}\sum_{b+c\le a}
\sum_{\begin{smallmatrix}b+c\le a\\ c\ne a\end{smallmatrix}}
\|\partial\Gamma^bu^i\partial^2\Gamma^cu^j\partial\Gamma^au^k\|_{L^1(r\le
c_0t)}.
\]
Since $r\le c_0 t$, we have that
 $\langle c_it-r\rangle\le
C \langle t \rangle$ for each $i=1,\ldots,m$.  Thus, using
\eqref{ellinf6}, a typical term
can be estimated by 
\begin{multline*}
C\langle t \rangle^{-3/2}
\|\langle c_it-r\rangle^{1/2} \partial\Gamma^bu^i
\langle c_jt-r\rangle\partial^2\Gamma^cu^j
\partial\Gamma^au^k\|_{L^1(r\le c_0t)}\\
\le C \langle t \rangle^{-3/2}
\|\langle c_it-r\rangle^{1/2}\partial\Gamma^bu^i\|_{L^\infty}
\|\langle c_jt-r\rangle\partial^2\Gamma^cu^j\|_{L^2}
\|\partial_t\Gamma^au^k\|_{L^2}\\
\le C \langle t \rangle^{-3/2}
\Big[E_{|b|+3}^{1/2}(u(t))+\mathcal{X}_{|b|+3}(u(t))\Big]
\mathcal{X}_{|c|+2}(u(t))E_\mu^{1/2}(u(t)).
\end{multline*}
In the preceding, we have $|b|+3\le\od$, $|c|+2\le\mu$, and $|a|+1\le\mu$.
With the aid of 
Lemma \ref{we}, we have achieved an upper bound of the form
\[
C \langle t \rangle^{-3/2} E_\mu(u(t))E_\od^{1/2}(u(t))
\]
for the portion of the integrals over $r\le c_0 t$
on the right of \eqref{enmain1}.

\subsubsection*{Away from the origin}
It remains to estimate the right-hand side of \eqref{enmain1}
for $r\ge c_0t$.  

First, we consider the nonresonant terms, i.e.\ those for
which $(i,j,k)\ne(k,k,k)$.  If $i\ne j$ and $r\ge c_0t$, then
$\langle t \rangle^{3/2}
\le C \langle r \rangle \langle c_it-r\rangle^{1/2}
\langle c_jt-r\rangle$.  
Using \eqref{ellinf6}
we have the estimate
\begin{multline*}
\|\partial\Gamma^bu^i
\partial^2\Gamma^cu^j
\partial\Gamma^au^k\|_{L^1(r\ge c_0t)}\\
\le C\langle t \rangle^{-3/2}
\|\langle r \rangle \langle c_it-r\rangle^{1/2}
\partial\Gamma^bu^i\|_{L^\infty}
\|\langle c_jt-r\rangle \partial^2\Gamma^cu^j\|_{L^2}
\|\partial\Gamma^au^k\|_{L^2}\\
\le C\langle t \rangle^{-3/2}
 \Big[E_{|b|+3}^{1/2}(u(t))+\mathcal{X}_{|b|+3}(u(t))\Big]
\mathcal{X}_{|c|+2}(u(t))
E_{|a|+1}^{1/2}(u(t))\\
\le C\langle t \rangle^{-3/2}
 E_\mu(u(t)) E_\od^{1/2}(u(t)).
\end{multline*}
Otherwise, if $j\ne k$, we pair the weight $\langle r \rangle \langle
c_kt-r\rangle^{1/2}$ with $\partial\Gamma^2u^k$ in $L^\infty$
to get the same upper bound.

We are left to consider the resonant terms in \eqref{enmain1}, 
i.e.\ $(i,j,k)= (k,k,k)$, 
in the region $r\ge c_0t$.  It is here, finally, where
the null condition enters.
An application of Lemma \ref{ptwnc} yields the following upper bound
for these terms:
\begin{multline*}
C\langle t\rangle^{-1}\sum_k
\sum_{\begin{smallmatrix}b+c=a\\c\ne a\end{smallmatrix}}
\Big[
\|\Gamma^{b+1}u^k \partial^2\Gamma^cu^k 
\partial\Gamma^au^k\|_{L^1(r\ge c_0t)}\\
+\|\partial\Gamma^{b}u^k \partial\Gamma^{c+1}u^k
\partial\Gamma^au^k\|_{L^1(r\ge c_0t)}\\
+\|\langle c_kt-r\rangle \partial\Gamma^{b}u^k \partial^2\Gamma^cu^k
\partial\Gamma^au^k\|_{L^1(r\ge c_0t)}\Big].
\end{multline*}
We still need to squeeze out an additional
decay factor of $\langle t \rangle^{-1/2}$.

Since $r\ge c_0t$, we have 
$\langle r \rangle \le C \langle t \rangle$.
Thus, we have using \eqref{ellinf4}
\begin{align*}
\|\Gamma^{b+1}u^k& \partial^2\Gamma^cu^k 
\partial\Gamma^2u^k\|_{L^1(r\ge c_0t)}\\
&\le C \langle t \rangle^{-1/2} 
\|\langle r \rangle^{1/2}\Gamma^{b+1}u^k\|_{L^\infty(r\ge c_0t)}
\|\partial^2\Gamma^cu^k\|_{L^2}
\|\partial\Gamma^au^k\|_{L^2}\\
&\le C \langle t \rangle^{-1/2} 
E_{|b|+3}^{1/2}(u(t)) E_\mu(u(t))\\
&\le C \langle t \rangle^{-1/2}
E_{\od}^{1/2}(u(t)) E_\mu(u(t)).
\end{align*}

In a similar fashion, the second term is handled using
\eqref{ellinf5}:
\begin{align*}
\|\partial&\Gamma^{b}u^k \partial\Gamma^{c+1}u^k
\partial\Gamma^2u^k\|_{L^1(r\ge c_0t)}\\
&\le C \langle t \rangle^{-1}
\|\partial\Gamma^{b}u^k\|_{L^2}
\|\langle r \rangle\partial\Gamma^{c+1}u^k\|_{L^\infty(r\ge c_0t)}
\|\partial\Gamma^au^k\|_{L^2}\\
&\le C \langle t \rangle^{-1}
E_{|c|+3}^{1/2}(u(t)) E_\mu(u(t))\\
&\le C \langle t \rangle^{-1}
E_{\od}^{1/2}(u(t)) E_\mu(u(t)).
\end{align*}

The final set of terms are estimated using \eqref{ellinf5}
again and \eqref{we1}.
\begin{align*}
\|\langle c_kt&-r\rangle \partial\Gamma^{b}u^k \partial^2\Gamma^cu^k
\partial\Gamma^au^k\|_{L^1(r\ge c_0t)}\\
&\le C \langle t \rangle^{-1}
\langle r \rangle \partial\Gamma^{b}u^k\|_{L^\infty(r\ge c_0t)}
\|\langle c_kt-r\rangle  \partial^2\Gamma^cu^k\|_{L^2}
\|\partial\Gamma^au^k\|_{L^2}\\
&\le C \langle t \rangle^{-1}
E_{|b|+3}^{1/2}(u(t))\mathcal{X}_{|c|+2}(u(t))E_\mu^{1/2}(u(t))\\
&\le C \langle t \rangle^{-1}E_{\od}^{1/2}(u(t)) E_\mu(u(t)).
\end{align*}

Combining all the estimates in this subsection, we obtain,
thanks to \eqref{enequiv},
the following inequality for the lower energy:
\begin{align}
\label{en3}
\widetilde{E}'_\mu(u(t))
&\le C \langle t \rangle^{-3/2} E_\mu(u(t)) E_\od^{1/2}(u(t))\\
\nonumber
&\le C \langle t \rangle^{-3/2} \widetilde{E}_\mu(u(t)) 
\widetilde{E}_\od^{1/2}(u(t)).
\end{align}

\subsection*{Conclusion of the proof}

By \eqref{enequiv}, we have that the modified energy
satisfies
$\widetilde{E}_\mu^{1/2}(u(t)) \le C\eps$ for $0\le t < T_0$.
So from \eqref{en2}, we find that
\[
\widetilde{E}_\od(u(t))\le \widetilde{E}_\od(u(0))
\langle t \rangle^{C\eps},
\]
provided $\eps$ is small.  Inserting this bound into
\eqref{en3} and using \eqref{enequiv}, we obtain
\begin{multline*}
(1/2)
{E}_\mu(u(t))\le
\widetilde{E}_\mu(u(t))\le \widetilde{E}_\mu(u(0))
\exp CI\widetilde{E}_\od^{1/2}(u(0))\\
\le 2{E}_\mu(u(t)) \exp 2CI{E}_\od^{1/2}(u(0))\le2\eps^2, 
\end{multline*}
with $I= \int_0^\infty \langle s \rangle^{-3/2+C\eps} ds$.
With this we see ${E}^{1/2}_\mu(u(t))$ remains strictly lees than
$2\eps$ throughout the closed interval $0\le t \le T_0$.
This shows that $E_\od(u(t))$ is bounded for all time, which
completes the proof of Theorem \ref{thm}.


\end{document}